\documentclass[11pt]{article}

\usepackage{amssymb}
\usepackage{amsmath}
\usepackage{amsfonts}

\newcommand{\F}{\mathbb{F}_q}
\newcommand{\Fd}{\mathbb{F}_{q^d}}
\newcommand{\Fk}{\mathbb{F}_{q^k}}
\newcommand{\Fn}{\mathbb{F}_{q^n}}
\newtheorem{thm1}{Theorem} 
\newtheorem{cor1}{Corollary}

\newtheorem{prop1}{Proposition}
\newtheorem{lem1}{Lemma}

\newenvironment{pf}{\begin{trivlist}\item[]{\bf Proof. }}%
{\samepage\hfill$\diamond$\end{trivlist}}

\begin{document}
\title{ Irreducible Compositions of Polynomials\\ over Finite Fields}

\author{Melsik K. Kyuregyan\footnote{Institute for Informatics and Automation Problems,
National Academy of Sciences of Armenia,
 P. Sevak street 1, Yerevan 0014, Armenia; \texttt{melsik@ipia.sci.am} }
 \and Gohar M. Kyureghyan\footnote{Department of Mathematics, Otto-von-Guericke University, Universit\"atsplatz 2, 39106 Magdeburg, Germany; \texttt{gohar.kyureghyan@ovgu.de}}}
\date{}
\maketitle{}

\begin{abstract}
This paper  is devoted to the composition method of constructing  families of irreducible polynomials over  finite fields.

\end{abstract}

\section{Introduction}\label{1Introduction}

Let $d$ be a divisor of $n$.
It is well known that an irreducible polynomial over $\F$ of degree $n$ splits into $d$ distinct irreducible factors of degree $n/d$ over $\Fd$.
Moreover, if $\displaystyle g(x) = \sum_{i=0}^{n/d}a_i x^i \in \Fd[x]$ is a factor of $f(x)$, then the remaining factors are 
$$
g^{(u)}(x)= \displaystyle \sum_{i=0}^{n/d} a_i^{q^u} x^i,$$ 
where $1 \leq u \leq d-1$. Consequently, the factorization of $f(x)$ in $\Fd[x]$ is given by
\begin{equation}\label{main}
f(x) = \prod_{u=0}^{d-1}g^{(u)}(x),
\end{equation}
where the notation $g(x) = g^{(0)}(x)$ is used. The converse of this statement is not true: Given an irreducible polynomial of degree $n/d$ over $\Fd$
the product $\displaystyle \prod_{u=0}^{d-1}g^{(u)}(x)$ is a polynomial over $\F$, but it must not necessarily be irreducible over $\F$.
To ensure that this product is irreducible over $\F$ it must be requested that
$\Fd$ is the smallest extension of $\F$ containing the coefficients of $g(x)$.
More precisely, it holds:

\begin{lem1}
\label{lemma3N}
A monic  polynomial $f(x) \in \mathbb{F}_q[x]$ of degree $n=dk$ is \mbox{ irreducible} over $\mathbb{F}_q$ if and only if there is a monic irreducible polynomial $g(x)= \sum_{i=0}^k g_u x^u$ over $\mathbb{F}_{q^d}$ of degree $k$ 
such that  $\F(g_0, \ldots, g_k) = \Fd$ and
 $f(x)= \displaystyle \prod_{v=0}^{d-1}g^{(v)} (x)$ in $\mathbb{F}_{q^d}[x]$.
\end{lem1}

As shown in Section \ref{2Preliminaries}, given an irreducible polynomial of degree $n$ over $\F$ and suitable elements in $\Fk$,
Lemma \ref{lemma3N} implies the following construction of irreducible polynomials of degree $nk$ over $\F$:

\begin{thm1}
\label{theorem1}
Let  $n>1$, $\gcd (n, k)=1$ and $f(x)$ be an irreducible polynomial of degree $n$ over   $\mathbb{F}_{q}$. Further, let $\alpha \ne 0$ and $\beta$ be elements of $\mathbb{F}_{q^k}$. 
Set $g(x) := f(\alpha x+ \beta)$. Then the polynomial 
\begin{equation}
\label{eq1}
F(x)=\prod _{a=0}^{k-1} g^{(a)} (x) 
\end{equation}
of degree $nk$
is irreducible over  $\mathbb{F}_{q}$ if and only if $\F(\alpha, \beta) = \Fk$.
\end{thm1}

The problem of reducibility of polynomials over finite fields is a case of special interest and plays an important role in modern engineering \cite{5Albert,7Cohen,1Lidl,8Menezes,7Zierler}. 
One of the  methods for constructing irreducible polynomials  is the composition method which allows constructions of irreducible polynomials of higher degree from the given  
irreducible polynomials with the use of  a substitution operator  (see \cite{2Cohen,2Kyuregyan,Scheerhorn}). Probably the most powerful result in this area is the following theorem by S. Cohen:

\begin{thm1}[Cohen~\cite{1Cohen}]\label{cohen}
Let $f(x), g(x) \in  \mathbb{F}_q[x]$  be relatively prime polynomials and let $P(x)  \in  \mathbb{F}_q[x]$ be an irreducible polynomial of degree $n$. Then the composition 
$$
F(x)=g^n(x)P\big( f(x)/g(x) \big) 
$$ 
is irreducible over $\mathbb{F}_q$
 if and only if $f(x)-\alpha g(x)$  is irreducible over $\mathbb{F}_{q^n}$ for a zero $\alpha \in \mathbb{F}_{q^n}$ of $P(x)$.
\end{thm1}

Theorem \ref{cohen} was employed by  several authors, including Chapman \cite {Chapman}, Cohen \cite {2Cohen},  
McNay \cite{McNay}, Meyn \cite{Meyn}, Scheerhorn \cite{Scheerhorn} and Kyuregyan \cite{1Kyuregyan}--\cite {3Kyuregyan} to give iterative constructions of irreducible polynomials and N-polynomials over finite fields.
Observe that Lemma \ref{lemma3N} yields a proof for Theorem \ref{cohen}. Indeed, over $\Fn$ the polynomial $P(x)$ is the product $\displaystyle \prod_{i=0}^{n-1}(x-\alpha ^{q^i})$ and thus
$$
F(x) = g^n(x)P\big( f(x)/g(x) \big) = \prod_{i=0}^{n-1}\left(f(x) - \alpha ^{q^i}\,g(x)\right) =  \prod_{i=0}^{n-1}\left(f(x) - \alpha \,g(x)\right)^{(i)}.
$$

In Section \ref{3Composition}
  we apply Theorem \ref{theorem1} to construct explicit families of  irreducible polynomials over finite fields.

In particular, using the results  by Ore-Gleason-Marsh \cite {7Zierler}, Dickson \cite{5Albert},  Sidelnikov \cite{6Sidelnikov}  we obtain  explicit  families of  irreducible polynomials of degrees $n(q^m-1)$ 
and $n(q^n+1)$ over $\mathbb{F}_q$  from a given irreducible polynomial of degree $n$  and  a primitive polynomial of degree $m$ over $\F$.  


\section{Preliminaries}
\label{2Preliminaries}

Throughout this paper we assume, without loss of generality, that the considered polynomials  are monic, i.e. with the leading coefficient $1$.
Let $f(x)$ be a monic irreducible polynomial of degree $n$ over $\mathbb{F}_q$
and let $\beta$ be a zero of $f(x)$. The field $\mathbb{F}_q (\beta)=\mathbb{F}_{q^n}$ is an $n$-dimensional extension of $\mathbb{F}_q$, which is a vector space of dimension $n$ over $\mathbb{F}_q$.

We say that the degree of an element $\alpha $ over $\mathbb{F}_q$ is equal to $k$ and write  $\deg_q(\alpha)=k$ if
$\mathbb{F}_q(\alpha)$ is a $k$-dimensional vector space over $\mathbb{F}_q$.
An element  $\alpha \in \mathbb{F}_{q^k}$   is called a proper element of $\mathbb{F}_{q^k}$ over $\mathbb{F}_{q}$ if  $\deg_q(\alpha)=k$, which is equivalent to the property that $\alpha \not \in \mathbb{F}_{q^v}$ for any proper divisor $v$ of $k$.
Similarly,  we say that the degree of a subset $A=\{ \alpha_1, \alpha_2, \cdots, \alpha_r\} \subset \mathbb{F}_{q^k}$ over $\mathbb{F}_q$ is equal to $k$
 and  write $\deg_q(\alpha_1, \alpha_2, \cdots, \alpha_r)=k$,
 if for any proper divisor $v$ of $k$ there exists at least one element $\alpha_u \in A$ such that $\alpha_u \not \in \mathbb{F}_{q^v}$.\footnote{A proper divisor of a natural number $n$ is a divisor of  $n$ other than $n$ itself.
}


The following results are well known and can be found for example in  \cite{1Lidl}. 
   
\begin{prop1}[\cite{1Lidl}, Theorem 3.46]
\label{Statement1N} 
Let $f(x)$ be a monic irreducible polynomial of degree $n$ over $\mathbb{F}_q$ and let  $k \in N$. Then $f(x)$ factors into $d$ irreducible polynomials in $\mathbb{F}_{q^k}[x]$ of the same degree $nd^{-1}$, where $d=gcd(n,k)$.
 \end{prop1}

\begin{prop1}[\cite{1Lidl}, Corollary 3.47]
\label{Statement2N} 
An irreducible polynomial over $\mathbb{F}_q$ of degree $n$ remains irreducible over extension field  $\mathbb{F}_{q^k}$  of  $\mathbb{F}_q$ if and only if $n$ and $k$ are relatively prime. 
 \end{prop1}

\begin{prop1}[\cite{1Lidl}, Theorem 3.29]
\label{Statement3N} 
The product $I(q,n;x)$ of all monic irreducible polynomials of degree $n$ in $\mathbb{F}_{q}[x]$ is given by 
$$
I(q,n;x)=\prod_{d\mid n}(x^{q^d}-x)^{\mu(n/d)}=\prod_{d\mid n}(x^{q^{n/d}}-x)^{\mu(d)},
$$
where $\mu (x)$ is  the M\"{o}ebius function.
 \end{prop1}

Given $0 \leq a\leq k-1$ and $g(x)=\displaystyle \sum_{u=0}^{m} b_u x^u \in \mathbb{F}_{q^k}[x]$, we use the notation
 $$g^{(a)}(x)= \displaystyle \sum_{u=0}^{m} b_u^{q^a} x^u.$$

The following lemma is well known and is an immediate consequence of Proposition \ref{Statement1N}. 

\begin{lem1}
\label{lemma1N}
Let   $f(x)$ be a monic irreducible polynomial of degree $dk$  over $\mathbb{F}_{q}$. Then there is 
   a monic irreducible divisor $g(x)$ of degree $k$ of  $f(x)$ in $\mathbb{F}_{q^d}[x]$. Moreover, every irreducible factor of $f(x)$ in $\mathbb{F}_{q^d}[x]$ is given  by $g^{(v)}(x)$ for some $0\leq v \leq d-1$.
In particular, the factorization of $f(x)$ in $\Fd[x]$ is 
\begin{equation}
\label{eq1N}
f(x)=\prod_{v=0}^{d-1}g^{(v)}(x).
\end{equation}
\end{lem1}

It is easy to see that, in general, the converse of Lemma \ref{lemma1N} does not hold. To ensure the converse statement, a factor $g(x)$ must be described more precisely, as
it is done in Lemma \ref{lemma3N} stated in Introduction.

{\bf PROOF of Lemma \ref{lemma3N}.}
Suppose $f(x)$ is irreducible over $\mathbb{F}_q$. Then 
 by Lemma \ref{lemma1N}\  there is an irreducible polynomial $g(x) = \displaystyle \sum_{u=0}^k g_u x^u $ of degree $k$ over $\mathbb{F}_{q^d}$  such that
\begin{equation}
\label{eq3N}
f(x)=\prod_{v=0}^{d-1}g^{(v)} (x)
\end{equation}
over $\mathbb{F}_{q^d}$.
Next we show that  the degree of  the set of  coefficients of $g(x)$ over $\mathbb{F}_{q}$  is equal to $d$. 
Suppose, on the contrary that $\deg _q ( g_0, g _1, \ldots, g_k )=s$, where $d= rs$ and $s < d$. 
Then, because of $\mathbb{F}_{q^s}[x] \subset \mathbb{F}_{q^d}[x]$,the polynomial $g(x)$ is  also  irreducible over $\mathbb{F}_{q^s}$ and by 
  Lemma \ref{lemma1N}  
\begin{equation}
\label{eq4N}
f(x)=\prod_{w=0}^{s-1}h^{(w)} (x)
\end{equation}
over $\mathbb{F}_{q^s}$ and 
$h^{(w)}(x) = \displaystyle \sum_{u=0}^{rk}h_u^{q^w}x^u $, $w=0,1,2,\ldots, s-1$, are distinct irreducible polynomials of degree $rk$ over  $\mathbb{F}_{q^s}$.
Combining (\ref{eq3N}) and (\ref{eq4N}) we get
$$
f(x) = \prod_{w=0}^{s-1} h^{(w)} (x) = \prod_{v=0}^{d-1} g^{(v)}(x)
$$
in $\mathbb{F}_{q^d}[x]$, which  contradicts to the uniqueness of the decomposition into irreducible factors in $\mathbb{F}_{q^d}[x]$.

To prove the converse, let  $g(x)$ be an  irreducible polynomial of degree $k$ over $\mathbb{F}_{q^d}$ and let $\alpha \in \mathbb{F}_{q^{dk}}$ be a zero of $g(x)$.  By Proposition  \ref{Statement3N}  
$$
I(q,dk;x)= \Big( x^{q^{dk}}-x \Big) \prod_{\stackrel {\scriptstyle \delta \mid dk}{\scriptstyle \delta \not =dk}} \Big( x^{q^{\delta}}-x \Big) ^{\mu (dk/ \delta)},
$$
which yields
$$
I(q,dk, \alpha)= \Big( \alpha^{q^{dk}}-\alpha\Big) \prod_{\stackrel {\scriptstyle \delta \mid dk}{\scriptstyle \delta \not =dk}} \Big( \alpha^{q^{\delta}}-\alpha \Big) ^{\mu (dk/ \delta)}=0,
$$
since $\alpha^{q^{dk}}= \alpha$.
Thus, $\alpha$ is a zero of  $I(q,dk, x)  \in \mathbb{F}_{q}[x]$ implying that $g(x)$ divides $ I(q,dk, x)$ in $\mathbb{F}_{q^d}[x]$.
In particular, there exists an irreducible polynomial $f(x)$ of degree $dk$  over $\mathbb{F}_{q}$ which is divisible by $g(x)$ in $\mathbb{F}_{q^d}[x]$. From Lemma \ref {lemma1N} it follows that $f(x)$ factors as
$$
f(x) = \prod_{v=0}^{d-1} g^{(v)}(x)
$$
in the ring $\mathbb{F}_{q^d}[x]$. 

Later we will use the following easy consequence of Proposition  \ref{Statement2N}.

\begin{lem1}
\label{lemma1}
Let  $\gcd(n, k)=1$,  \ $f(x)$ be an irreducible polynomial of degree $n$ over $\mathbb{F}_{q}$ and let $\alpha\ne 0, \beta \in \mathbb{F}_{q^k}$. Then the polynomial $g(x)=f(\alpha x+\beta)$ is irreducible over $\mathbb{F}_{q^k}$.
\end{lem1}

The next lemma provides the conditions on the elements $\alpha, \beta$ under which the degree of the set of coefficients  of  $g(x)=f(\alpha x+\beta)$ is equal to $k$ over $\F$.

\begin{lem1}
\label{lemma2}
Let   $n>1$ 
and $f(x)$ be an irreducible polynomial of degree $n$ over  \ $\mathbb{F}_{q}$.
Further, let \ $\gcd(n, k)=1$  and let   $\alpha, ~\beta \in \mathbb{F}_{q^k},~ \alpha \ne 0$. 
 Then the degree of the set of coefficients $\{g_0, g_1, \ldots, g_n\}$ of the polynomial $g(x)=f(\alpha x+\beta)$ is equal to $k$ over  $\mathbb{F}_{q}$ if and only if $\deg_q(\alpha, \beta)=k$.
\end{lem1}

\begin{pf} Suppose $\deg_q(\alpha, \beta)=k$.
Let $\theta \in \mathbb{F}_{q^n}$ be a zero of $f(x)$. 
Then $\gamma = \alpha_1\theta + \alpha_2  \in \mathbb{F}_{q^{nk}}$ is a zero of $g(x)$, where $\alpha_1=\alpha^{-1}$ and $\alpha_2=-\alpha^{-1}\beta$.
Suppose,  that  \ the   \ degree  \ of  \ the  \ set \ of  \ coefficients $\{g_0, g_1, \ldots, g_n\}$ of $g(x)$ is $v$ over $\mathbb{F}_{q}$, where $1 \leq v \leq k$ divides $k$. 
Hence  $\gamma$ is a root of the irreducible polynomial $g(x)$ of degree $n$ over $\mathbb{F}_{q^v}$, and therefore $\gamma$  a proper element of  $\mathbb{F}_{q^{nv}}$ over $\mathbb{F}_{q^{v}}$. In particular, it holds
\begin{equation}\label{eq-ld}
\gamma ^{q^{nv}} = (\alpha_1\theta + \alpha_2)^{q^{nv}} = \alpha_1 ^{q^t}\theta + \alpha_2^{q^t} = \gamma = \alpha_1\theta + \alpha_2,
\end{equation}
where $nv \equiv t \pmod k$ and $0 \leq t \leq k-1$. To prove the statement of the lemma, we must show that $t=0$. Suppose, to the contrary that $1 \leq t \leq k-1$.
From (\ref{eq-ld}) it follows that
$$
(\alpha_1^{q^t} - \alpha_1) \cdot \theta + (\alpha_2^{q^t}- \alpha_2)\cdot 1 =0.
$$
Since $\theta$ and $1$ are linearly independent over $\mathbb{F}_{q^k}$, the latter identity implies
$$
\alpha_1^{q^t} - \alpha_1 = 0 ~\textrm{ and }~ \alpha_2^{q^t} - \alpha_2 = 0.
$$
Hence $\alpha_1, \alpha_2 \in \mathbb{F}_{q^s}$ with $s = \gcd(k,t) < k$.
This yields that $\alpha \in  \mathbb{F}_{q^s}$ and $-\alpha \cdot \alpha_2 = \beta \in  \mathbb{F}_{q^s}$, and thus $ \mathbb{F}_{q}(\alpha, \beta) =  \mathbb{F}_{q^s}$, contradicting to the assumption
that  $ \mathbb{F}_{q}(\alpha, \beta) =  \mathbb{F}_{q^k}$.
\end{pf}

Observe that   Lemmas \ref{lemma3N}\,-\,\ref{lemma2}  imply the statement of Theorem \ref{theorem1} stated in the introduction.

\section{Irreducibility of Polynomial Compositions}
\label{3Composition}

In this section we apply Theorem \ref{theorem1} to describe several explicit families of irreducible polynomials  over  $\mathbb{F}_{q}$. 
We start by showing that Theorem \ref{theorem1} implies a proof for a result stated by Varshamov in \cite{3Varshamov} with no proof.

Recall that given $l,\,m$ with $\gcd(l,m)=1$,  the natural number $o\ne 0$ is called the order of $l$ modulo $m$ if it is the minimal number satisfying $l^o \equiv 1 \pmod m$. 

\begin{thm1}[Varshamov~\cite{3Varshamov}]
\label{corollary1}
Let $r$ be an odd  prime number which does not divide $q$  and  $r-1$ be the order of $q$ modulo $r$.
Further, let $n >1,~\gcd(n, r-1)=1$ and   $f(x)$  be an irreducible  polynomial of degree  $n$ over  $\mathbb{F}_{q}$ belonging to order $t$.
Define the polynomials $R(x)$ and $\psi(x)$ over $\mathbb{F}_q$ as follows:
Set  $x^r \equiv R(x) \! \! \pmod {f(x)}$ and $ \psi (x) =\displaystyle \sum_{u=0}^{n} \psi_u x^{u}$, where $\psi(x)$ is the nonzero polynomial of minimal degree satisfying the congruence
\begin{equation}
\label{definition-psi}
\displaystyle \sum_{u=0}^{n} \psi _u (R(x))^u \equiv 0 \!\!\! \pmod{f(x)}.
\end{equation} 
 Then the polynomial $\psi(x)$ is an irreducible polynomial of degree $n$ over $\mathbb{F}_{q}$ and 
$$
F(x) =f^{-1}(x) \: \psi (x^r)
$$
is an irreducible polynomial of degree $(r-1)n$  over $\mathbb{F}_q$. Moreover $F(x)$  belongs to order $rt$.
\end{thm1}

\begin{pf}
Let $\alpha \in \mathbb{F}_{q^n}$ be a zero of $f(x)$. Then  $x^r \equiv R(x) \! \! \pmod {f(x)}$ is equivalent to $\alpha ^r = R(\alpha)$ in $\mathbb{F}_{q^n}$.
Note that the condition that  $\psi(x)$ is the nonzero polynomial of minimal degree satisfying (\ref{definition-psi}) implies that
$\psi(x)$ is the minimal polynomial of $R(\alpha) = \alpha^r$ over $\mathbb{F}_q$. In particular, $\psi(x)$ is irreducible over $\mathbb{F}_q$.
In order to prove that the degree of $\psi$ is $n$, we will show that $\alpha ^r$ is a proper element of $\mathbb{F}_{q^n}$ over $\mathbb{F}_{q}$  by proving that
 the (multiplicative) order of $\alpha^r$ is equal to the one of $\alpha$.
By the assumption on $f(x)$ the  order of $\alpha$ is $t$. Thus
the order of $\alpha^r$ is $t/\gcd(t,r)$ and it is enough  to show that $\gcd(t,r)=1$. To prove the latter recall that the smallest $i$ such that $r$ 
divides $q^{i}-1$ is  $r-1 \ne 1$, further $t$ divides $q^n-1$ and finally 
$$
\gcd(q^n-1, q^{r-1}-1) = q^{\gcd(n,r-1)}-1 = q-1.
$$
Now we consider the polynomial $F(x) = \psi(x^r) f^{-1}(x)$. Over  $\mathbb{F}_{q^n}$ we have
$$
f(x) =\displaystyle \prod_{u=0}^{n-1} (x-\alpha^{q^u}) ~\mbox{  and }~
\psi(x)= \displaystyle \prod_{u=0}^{n-1}(x-\alpha^{rq^u})
$$
and consequently
$$
F(x)=\prod_{u=0}^{n-1}\frac{x^r-\alpha^{rq^u}}{x-\alpha^{q^u}}  =\prod_{u=0}^{n-1}\Big(x^{r-1}+\alpha^{q^u} x^{r-2}+\cdots+ \alpha^{q^u(r-2)}x+\alpha^{q^u(r-1)}\Big).
$$
Set 
$$
g(x) : =x^{r-1}+\alpha x^{r-2}+ \cdots +\alpha^{r-2}x+ \alpha^{r-1}.
$$
Then $F(x)  =\displaystyle \prod_{u=0}^{n-1} g^{(u)}(x)$. 
 Note that $g(x)=\alpha^{r-1}h (\alpha^{-1}x)$, where $h(x)=x^{r-1}+x^{r-2}+ \cdots + x+1$. 
It is well known that the polynomial $h(x)$ is irreducible over $\mathbb{F}_{q}$ if and only if $r$ is a prime number and the order of $q$ modulo $ r$ is $r-1$.
Hence the irreducibility of $F(x)$ over $\F$ is implied by Theorem \ref{theorem1}.

To complete the  proof it remains to show that the order of $F(x)$ is $rt$. 
Let $\beta$ be a zero of $h(x)$. Since $x^r-1=(x-1)h(x)$, the order of $\beta$ is $r$. From $F(x)  =\displaystyle \prod_{u=0}^{n-1} g^{(u)}(x)$ and $g(x) = \alpha^{r-1}h (\alpha^{-1}x)$
it follows that the element $\alpha \beta$ is a zero of $F(x)$. Now the statement follows from the fact that the order of $\alpha \beta$ is the smallest common multiple of the orders
of $\alpha$ and $\beta$, {\em i.e.} $rt$ since $\gcd(r,t)=1$ as shown above.
\end{pf}

Recall that a polynomial $\displaystyle l(x) = \sum_{i=0}^n a_ix^{q^i} \in \mathbb{F}_{q}[x]$ is called a linearized polynomial over $\mathbb{F}_q$. The polynomials 
$$
l(x) = \sum_{i=0}^{n} a_i x^{q^i}  ~\textrm{ and }~  \bar{l}(x) = \sum_{i=0}^{n} a_i x^{i} 
$$
are called $q$-associates of each other. More precisely, $\bar{l}(x)$ is the conventional $q$-associate of ${l}(x)$, and ${l}(x)$ is the linearized $q$-associate of $\bar{l}(x)$.

\begin{thm1}[Ore-Gleason-Marsh,  \cite {7Zierler}] 
\label{theorem2}
Let $f(x)=\displaystyle \sum_{u=0}^{n} a_u x^{u} \in \mathbb{F}_q[x]$ and $F(x)$ be its linearized $q$-associate. Then
the polynomial $f(x)$  is a primitive polynomial over $\mathbb{F}_q$ if and only if  the polynomial $x^{-1}F(x)=\displaystyle \sum_{u=0}^{n} a_u x^{q^u-1}$ is irreducible  over  $\mathbb{F}_q$.
\end{thm1} 

Given an irreducible polynomial of degree $n$ and a primitive polynomial of degree $m$ over $\F$, the next theorem yields an irreducible polynomial of degree $n(q^m-1)$ over $\F$.

\begin{thm1}
\label{theorem3}
Let $\gcd(n,q^m-1)=1$ and \ $l(x)=\displaystyle \sum_{v=0}^m b_v x^{q^v}$ such that its conventional $q$-associate $\bar{l}(x)\not =x-1$ is a primitive polynomial of degree $m$ over $\mathbb{F}_q$.
 Further, let $f(x)$  be an irreducible polynomial of degree $n$ over  $\mathbb{F}_{q}$. Define $R(x)$ and $\psi(x)$ as follows:
 $l(x) \equiv R(x) \! \! \pmod {f(x)}$ and $ \psi (x) =\displaystyle \sum_{u=0}^{n} \psi_u x^{u} \in \mathbb{F}_q[x]$ to be the nonzero polynomial of minimal degree satisfying the congruence 
\begin{equation}
\label{eq8N}
\sum _{u=0}^{n} \psi_u (R(x))^u \equiv 0 \!\!\! \pmod{f(x)}.
\end{equation}
Then $\psi(x)$  is an irreducible polynomial of degree $n$ over $\mathbb{F}_q$ and $F(x)=(f(x))^{-1}\psi(l(x))$ is an irreducible polynomial of degree $n(q^m-1)$ over $\mathbb{F}_q$.
\end{thm1}

\begin{pf}
First consider the case $n=1$, {\it{i.e.}} $f(x)=x+a$ with $a \in \mathbb{F}_{q}$.
Then
\begin{eqnarray*}
l(x) && =  x^{q^m} + b_{m-1}x^{q^{m-1}} + \cdots+ b_1 x^q + b_0 x   \\
& &  =(x+a)^{q^m}+b_{m-1}(x+a)^{q^{m-1}}+ \cdots + b_1(x+a)^q+b_0 (x+a)  \\
& &  -a (1+b_{m-1}+ \cdots + b_1+b_0) ,
\end{eqnarray*}
and, in particular,
$$
l(x) \equiv -a (1+b_{m-1}+ \cdots + b_1 +b _ 0) \! \!\! \pmod {(x+a)}. 
$$
Using the definition of $\psi(x)$ we get $\psi (x) =x+a(1+b_{m-1}+ \cdots + b_1 +b _ 0)$. And so  
\begin{eqnarray*}
F(x) & & = (f(x))^{-1} \psi(l(x))  \\
&& = \frac { x^{q^m} + b_{m-1}x^{q^{m-1}} + \cdots+ b_1 x^q + b_0 x + a (1+b_{m-1}+ \cdots + b_1+b_0)} {x+a}\\
&& = \frac { (x+a)^{q^m} + b_{m-1}(x+a)^{q^{m-1}} + \cdots + b_1 (x+a)^q + b_0 (x+a)} {x+a}  \\
&& = (x+a)^{q^m-1} + b_{m-1}(x+a)^{q^{m-1}-1} + \cdots + b_1 (x+a)^{q-1} + b_0.
 \end{eqnarray*}
The latter polynomial is irreducible over $\mathbb{F}_{q}$ by Theorem \ref {theorem2}.
 
We next consider the case $n>1$. 
Let $\alpha \in \Fn$ be a zero of $f(x)$.
 Consider the polynomial
$$
H(x)=x^{-1}l(x) =x^{q^m-1}+b_{m-1}x^{q^{m-1}-1}+ \cdots + b_1 x^{q-1}+b_0
$$
which is irreducible over $\mathbb{F}_{q}$ by Theorem \ref{theorem2}. Set  $h(x)= H(x-\alpha)$.
It is easy to see, that $h^{(u)}(x)= H(x-\alpha^{q^u})$ for $0\leq u \leq n-1$.
Using Theorem \ref{theorem1} we get that the polynomial
$$
F(x) = \prod_{u=0}^{n-1} h^{(u)}(x) = \prod_{u=0}^{n-1}H(x - \alpha ^{q^u})
$$
is irreducible over $\mathbb{F}_q$.

Note that by definition of $R(x)$ it holds $l(\alpha)=R(\alpha)$ in $\mathbb{F}_{q^n}$. Further, we have
\begin{eqnarray*}
f(x)F(x) & = & \prod_{u=0}^{n-1}(x-\alpha ^{q^u}) H(x - \alpha ^{q^u}) 
 = \prod_{u=0}^{n-1}(x-\alpha ^{q^u}) \frac{l(x - \alpha ^{q^u})}{x-\alpha ^{q^u}}\\
 & = &  \prod_{u=0}^{n-1}\left(l(x) - l(\alpha)^{q^u}\right) = \prod_{u=0}^{n-1}\left(l(x) - R(\alpha)^{q^u}\right).
\end{eqnarray*}

Observe that $\psi(x)$ is the minimal polynomial of $R(\alpha)$ over $\mathbb{F}_q$. Hence $\psi(x)$ is irreducible over $\mathbb{F}_q$.
It has degree $n$, since $R(\alpha)$ is a proper element of $\mathbb{F}_{q^n}$ over $\mathbb{F}_{q}$. 
Indeed, suppose on the contrary,
that the  degree of  $R(\alpha)$ over $\mathbb{F}_{q}$ is equal to $d$, where $d$ is a proper divisor of $n$. Then
$$
\prod_{u=0}^{n-1}(x-(R(\alpha))^{q^{u}})=\Big( \prod_{u=0}^{d-1} \left(x-\left(R(\alpha)\right)^{q^u} \right)\Big)^k  = (\psi(x))^k, 
$$
 where  $n=dk$.
Substituting $l(x)$ for $x$ in the expression above, we obtain
\begin{equation}\label{eq-contr}
f(x)F(x) =  \prod_{u=0}^{n-1} \Big( l(x) - (R(\alpha))^{q^u}\Big ) = \Big( \psi\big(l(x)\big) \Big)^k.
\end{equation}
Recall that  $f(x)$ and $F(x)$ are irreducible polynomials of degree $n$ and $n(q^m-1)$, resp., over  $\mathbb{F}_{q}$. Hence  (\ref{eq-contr}) forces  that 
$k=2$, $dq^m=n$ and $dq^m=n(q^m-1)$. In particular, it must hold $n=n(q^m-1)$, which is impossible, since  by assumption $\bar{l}(x) \not = x-1$, and therefore $q^m \ne 2$  and   $n(q^m-1) >n$. 
 
Finally it remains to note  that  (\ref{eq-contr}) holds with $k=1$, showing that 
$
F(x)= (f(x))^{-1}\psi( l(x)).
$
\end{pf}

Observe that the computing of the minimal polynomial $\psi(x)$ of $R(\alpha)$ in (\ref{eq8N}) is equivalent to solving a system of $n$ linear equations with $n$ unknowns $\psi_1, \ldots ,\psi_{n-1}$.

\vspace*{0.3cm}
For the choice $l(x) = x^q - \theta x$  Theorem  \ref{theorem3} yields:
\begin{cor1}
\label{theorem6}
Let $q>2$, \ $\gcd(n, q-1)=1$ and  $f(x)$ be an irreducible polynomial of degree $n$ over $\mathbb{F}_q$.
Further, let $\theta$ be a primitive element of  \ $\mathbb{F}_q$. Define $R(x)$ and $\psi(x)$ as follows: Let
 $x^{q}-\theta x \equiv R(x)\!\! \pmod {f(x)}$ and $\psi(x) =\displaystyle \sum_{u=0}^{n} \psi_ux^u$ to be the nonzero polynomial of  the least degree satisfying the congruence
\begin{equation}
\label{eq12N}
\sum_{u=0}^{n} \psi_u (R(x))^u \equiv 0 \!\! \pmod {f(x)}.
\end{equation}
Then  $\psi(x)$ is an irreducible polynomial of degree $n$ over $\F$ and $F(x)= \big(f(x)\big)^{-1} \psi(x^q-\theta x)$ is an irreducible polynomial of degree $n(q-1)$ over $\mathbb{F}_q$. 
\end{cor1}

Another consequence of Theorem \ref{theorem3} is:
\begin{cor1}
\label{cor1}
Let $\gcd(n, q^m-1)=1$, $l(x)=\displaystyle \sum_{v=0}^m b_vx^{q^v}$ such that
its convensional $q$-associate $\bar{l}(x) \ne x-1$ is  a primitive polynomial of degree $m$ over $\mathbb{F}_{q}$ and let $f(x)$ be an irreducible polynomial of degree $n$ over $\mathbb{F}_{q}$. 
For any $0 \leq i \leq n-1$ define
 $c_i=\displaystyle \sum_{u=0}^{\lfloor n^{-1}(m+1)\rfloor} b_{i+nu}$,
where $b_u = 0$ for $u >m$. Suppose there is an $i$
such that $c_i \ne 0$ and  $c_j=0$ for  $j\ne i, 0 \leq j \leq n-1$. Then the polynomial of degree $n(q^m-1)$
$$
F(x)=\big(f(x)\big)^{-1}f \big(c_i^{-1}l(x)\big)
$$
is irreducible over $\mathbb{F}_{q}$.
\end{cor1}

\vspace{-1cm}
\begin{pf}
We use the notation of Theorem \ref{theorem3}.
Clearly, we have $l(x) =\displaystyle \sum_{v=0}^m b_vx^{q^v}=\displaystyle \sum_{u=0}^{\lfloor n^{-1}(m+1)\rfloor} b_{i+nu}x^{q^{i+nu}}$.  
Let $\alpha \in \Fn$ be a zero of   $f(x)$. Then using the conditions on $c_i$ we get   
$R(\alpha)=\displaystyle \sum_{u=0}^{\lfloor n^{-1}(m+1)\rfloor} b_{i+nu}\alpha^{q^{i+nu}}=c_i\alpha^{q^i}$, implying that $\psi(x)=f(c_i^{-1}x)$.
Theorem \ref{theorem6} completes the proof.
\end{pf}

Next two examples  are applications of Corollary \ref{cor1}.

\vspace*{0.3cm}
\noindent
{\bf Example.}

\begin{description}
\item[(a)] Let $q=2$ and $n = 2$. Recall that the unique irreducible polynomail of degree $2$ over $\mathbb{F}_{2}$ is $f(x) = x^2 + x + 1$. 
Let $\bar{l}(x) = \sum_{v=0}^m b_vx^v$ be a primitive polynomial of degree $m$ over $\mathbb{F}_{2}$ and $l(x)$ its linearized $2$-associate. Then exactly one of the sums $c_0 = \sum_{j = 0}^{\lfloor{m+1}/2\rfloor}b_{2j}$
or $c_1 = \sum_{j = 0}^{\lfloor{(m+1)}/2\rfloor}b_{2j +1}$ is $0$, since $c_0 + c_1 = \bar{l}(1) = 1$. Hence by Corollary \ref{cor1} the polynomial
$$
\frac{l(x)^2 + l(x) + 1}{x^2 + x +1}
$$
is irreducible polynomial of degree $2(2^m-1)$ over  $\mathbb{F}_{2}$.

\item[(b)]
 Let $q=2$, $m=5$, $n=3$. The polynomial $\bar{l}(x)=x^5+x^4+x^2+x+1$ is primitive  over $\mathbb{F}_{2}$
 and the polynomial $f(x)=x^3+x+1$  is irreducible  over  $\mathbb{F}_{2}$. First, we compute $c_i$ from $\bar{l}(x)=\sum_{i=0}^{m}b_i x^i =x^5+x^4+x^2+x+1$:
$$
\begin{array}{c}
 c_0  =   b_0+b_3 = 1+0=1,  \\
 c_1  =  b_1+b_4 = 1+1=0, \\
 c_2  =  b_2+b_5 = 1+1=0.
\end{array}
$$
Hence, the assumptions of Corollary \ref{cor1} are fulfilled and thus  the polynomial $F(x)=(x^3+x+1)^{-1} \big( (l(x))^3+l(x)+1 \big)$, where $l(x)=x^{32}+x^{16}+x^4+x^2+x$,
or, more precisely, 
\begin {eqnarray*}
F(x)& = &\frac{(x^{32}+x^{16}+x^4+x^2+x)^3 +x^{32}+x^{16}+x^4+x^2+x+1}{x^{3}+x+1}=\\[2mm]
 &  & x^{93}+x^{91}+ x^{90}+ x^{89}+x^{86}+x^{84}+x^{83}+x^{82}+x^{79}+x^{77}+x^{76}+\\[4mm]
&  & x^{75} +x^{72}+x^{70}+ x^{69}+x^{68}+x^{65}+ x^{63}+x^{62}+x^{61}+x^{58}+x^{56}+ \\[4mm]
&  & x^{55}+x^{54}+x^{51}+x^{49}+ x^{48}+x^{47}+x^{45}+ x^{44}+x^{43}+x^{40}+x^{38}+\\[4mm]
 & & x^{37}+x^{36}+x^{33}+x^{31} + x^{30}+x^{27}+ x^{25}+x^{24}+x^{23}+ x^{20}+x^{18}+\\[4mm]
& & x^{17}+x^{16}+x^{9}+x^{7}+x^{6}+x^{5}+ x^{3}+x^{2}+1 
\end{eqnarray*}
is irreducible over  $\mathbb{F}_{2}$. 
\end{description}

Further we describe another  composition method that enables explicit constructions of  irreducible polynomials of degree $n(q^n-1)$  from a given primitive polynomial of degree $n$ over $\mathbb{F}_{q}$ by using a simple transformation. 
The method is based upon the following result.

\begin{thm1}[\cite{5Albert} Chapter V, Theorem 24 (Dickson's theorem)]
\label{theorem4}
Let $\theta$ be a primitive element of  $\mathbb{F}_{q}$,  $\beta$ be any  element of  $\mathbb{F}_{q}$, and $p^m>2$, where $m$ divides $s \; (q=p^s)$. Then the polynomial
$$
f(x)=x^{p^m}-\theta x+\beta
$$
is  the product of a linear polynomial and an irreducible polynomial of degree $p^m-1$ over $\mathbb{F}_{q}$.   
\end{thm1}

\begin{thm1}
\label{theorem5}
Let $q^n>2$, $\beta, \gamma \in \mathbb{F}_{q}$, $\beta  \not = - \gamma$ and  $f(x) \not = x-1$ be a primitive polynomial of degree $n$ over $\mathbb{F}_{q}$. Set $h(x)=f \big( (\beta+\gamma)x+1\big)$ and  $h^*(x)=x^n h\left(\frac{1}{x}\right)$. 
Then the polynomial 
$$
F(x)=(x-\gamma)^n f \big((x-\gamma)^{-1} (x^{q^n}+\beta)\big) \big(h^*(x-\gamma)\big)^{-1}
$$
is an irreducible polynomial of degree $n(q^n-1)$  over $\mathbb{F}_{q}$.   
\end{thm1}   

\begin{pf}
Let $\alpha \in \mathbb{F}_{q^n}$ be a zero of $f(x)$. Then in 
$\mathbb{F}_{q^n}[x]$ it holds  
\begin{equation}
\label{eq5}
f(x)=\prod_{u=0}^{n-1} \Big(x-\alpha^{q^u}\Big).
\end{equation} 
Substituting $(x-\gamma)^{-1} (x^{q^n}+\beta)$ for $x$ in (\ref{eq5}), and multiplying both sides of the equation by $(x-\gamma)^n$, we get
\begin{equation}
\label{eq6}
\big(x-\gamma \big)^n f \big((x-\gamma)^{-1}(x^{q^n}+\beta)\big)=\prod_{u=0}^{n-1}\Big(x^{q^n}-\alpha^{q^u}x +\beta +\gamma \alpha^{q^u}\Big).
\end{equation}
Since  $q^n>2$ and  $\alpha^{q^u}$ is a primitive element in $\mathbb{F}_{q^n}$,  Dickson's theorem yields that each of the polynomials $g^{(u)} = x^{q^n}-\alpha^{q^u}x+\beta+\gamma \alpha^{q^u}$ 
is product of a linear polynomial and an irreducible polynomial of degree $q^n-1$ over $\mathbb{F}_{q^n}$. 
Moreover, the linear factor of $g^{(u)}$ is $x-\theta^{q^u}$, where  $\theta^{q^u}=(\beta+\gamma \alpha^{q^u})(\alpha^{q^u}-1)^{-1}$, since $\theta^{q^u}$ is a zero of it.
Thus
$$
Q^{(u)}(x) = \frac{x^{q^n}-\alpha ^{q^u}x+\beta +\gamma \alpha^{q^u}}{x-\theta^{q^u}} = \frac{x^{q^n}- \theta^{q^{n+u}}-\alpha^{q^u} (x-\theta^{q^u})}{x-\theta^{q^u}}
$$
is irreducible over $\mathbb{F}_{q^n}$. Note that the free term of $Q^{(u)}(x)$ is $1-\alpha ^{q^u}$, and  in particular the degree of the set of its coefficients is $n$ over $\mathbb{F}_q$.
Consequently, by Lemma \ref{lemma3N} the polynomial $\displaystyle \prod_{u=0}^{n-1} Q^{(u)}(x)$ is irreducible over $\mathbb{F}_q$.
To complete the proof observe that
$$
F(x) = \frac{ (x-\gamma)^n f \big((x-\gamma)^{-1} (x^{q^n}+\beta)\big)  }{ \prod_{u=0}^{n-1} \big(x-\theta^{q^u}\big)} = \prod_{u=0}^{n-1} Q^{(u)}(x),
$$
since 
$$
\prod_{u=0}^{n-1} \big(x-\theta^{q^u}\big)=h^{*}(x-\gamma).
$$ 
Indeed, $\theta =(\beta+\gamma)(\alpha -1)^{-1}+\gamma$ and $(\beta+\gamma)^{-1}(\alpha -1)$ is a zero of  $h(x)=f\big((\beta+\gamma)x+1\big)$, which implies that $\theta$ is a zero of $h^*(x-\gamma)$.
\end{pf}

Further we obtain  explicit families of  irreducible polynomials of degree \mbox{$n(q^n+1)$} over finite fields using the following result: 
\begin{thm1}[Sidelnikov \cite{6Sidelnikov}] 
\label{theorem7}
Let  $w \in \mathbb{F}_q$ and $x_0 \in \mathbb{F}_{q^2} \setminus \mathbb{F}_q$ such that $x_0^{q+1}=1$. Then the   polynomial
$$
f(x) \: = \: x^{q+1} -wx^q -(x_0+x_0^q-w)x+1 \in \mathbb{F}_q[x]
$$
  is  irreducible  if and only if \  $\displaystyle \frac{w-x_0^q}{w-x_0}$ is a generating element of the multiplicative subgroup  $\Pi := \{ y \in \mathbb{F}_{q^2} ~|~ y^{q+1}=1 \}$ of $\mathbb{F}_{q^2}$. 
 Moreover, the polynomial $f(x)$ has linearly independent roots over $\mathbb{F}_q$.
\end{thm1}

\begin{thm1}
\label{theorem8}
Let $f(x)$ be an irreducible polynomial of degree $2n$ over $\mathbb{F}_q$ of order $e(q^n+1)$.
\begin{description}
\item[(a)] 
Let $\alpha \in \mathbb{F}_{q^{2n}}$ be a zero of $f(x)$. Set $\beta = \alpha ^{e}$.
Then the polynomial $x^{q^n+1}+x^{q^n}-(\beta^{q^n}+\beta +1) x+1$ is an irreducible polynomial over $\mathbb{F}_{q^n}$.
\item[(b)]
\vspace{0.2cm} 
Define the polynomials $R(x)$ and $\psi(x)$ over $\mathbb{F}_{q}$ as follows: Let
 $x^{eq^n} + x^e+1 \equiv R(x) \: (mod \: f(x))$ and $\psi(x)=\displaystyle \sum_{u=0}^{n} \psi _u x ^u$ be the nonzero polynomial of  the least degree satisfying the congruence  
\begin{equation}
\label{eq7}
\sum _{u=0}^n \psi_u (R(x))^u \equiv 0 \: (mod \: f(x)).
\end{equation}
Then the polynomial $\psi(x)$  is an irreducible polynomial of degree $n$  over $\mathbb{F}_q$. 
\item[(c)]
\vspace{0.2cm}
The polynomial $ F(x)=x^n\psi \left (\displaystyle \frac{x^{q^n+1}+x^{q^n}+1}{x} \right)$ is an irreducible polynomial  of degree $n(q^n+1)$ over $\mathbb{F}_q$.
\end{description}
\end{thm1}

\begin{pf}
{\bf (a)}  Note, that the order of $\beta$ is $q^n+1$, which does not divide $q^k-1$ for $k \leq n$. Hence  $\beta$ is a proper element of $\mathbb{F}_{q^{2n}}$ over $\mathbb{F}_{q}$.
Clearly $\gamma := \beta^{q^n}+\beta +1$ belongs to $\mathbb{F}_{q^n}$. Next we show that $\gamma$ is a proper element of $\mathbb{F}_{q^n}$ over $\mathbb{F}_{q}$.
Indeed, suppose $\gamma \in \mathbb{F}_{q^d}$ for some divisor $d$ of $n$. We have
$$
\gamma \beta = \beta^{q^n+1}+\beta^2 +\beta = 1 + \beta ^2 + \beta,
$$
and consequently, $\beta^2 +(1-\gamma)\beta +1=0$. Hence $\beta$ is a root of a quadratic polynomial over $\mathbb{F}_{q^d}$, implying that $[\mathbb{F}_{q^{2n}} : \mathbb{F}_{q^d}] \leq 2$
and thus $d=n$. To complete the proof of the statement (a), we show that the conditions of Theorem \ref{theorem7} are fullfiled. Indeed, choose $x_0 = \beta$ and $\omega = -1$.
It remains to note that $ \displaystyle \frac{\omega-x_0^{q^n}}{\omega-x_0}=\frac{-1-\beta^{q^n}}{-1-\beta}=\beta^{q^n} $ generates $\Pi$.\\
{\bf (b)} The congruence $x^{e{q^n}}+x^e+1 \equiv R(x) \!\! \pmod {(f(x))}$ is equivalent to the relation  $\alpha^{e{q^n}}+\alpha^e+1=R(\alpha) $ in  $\mathbb{F}_{q^{2n}}$ or $\beta^{q^n}+\beta +1 =R(\alpha)$.
Further, the condition that $\psi(x)$ is the nonzero polynomial of  the least degree satisfying  congruence (\ref{eq7}) is equivalent to the one that $\psi(x)$ is the minimal polynomial of $R(\alpha) = \beta^{q^n}+\beta +1$.
To complete the proof observe that the degree of $\psi(x)$ is $n$, since $\beta^{q^n}+\beta +1$ is a proper element of $\mathbb{F}_{q^n}$ over  $\mathbb{F}_{q}$ as shown in the proof of (a).\\
{\bf (c)} The polynomial $\psi(x)$ is the minimal polynomial of $ \beta^{q^n}+\beta +1$ over $\mathbb{F}_q$, and hence
\begin{equation}
\label{eq8}
\psi (x)=\prod _{u=0}^{n-1}(x-(\beta^{q^n}+\beta+1)^{q^u}).
\end{equation}
Substituting $\displaystyle \frac{x^{q^n+1}+x^{q^n}+1}{x}$ for $x$ in (\ref{eq8}), and multiplying both sides of the expression by $x^n$, we obtain 
$$
x^n \psi \left(\frac{x^{q^n+1}+x^{q^n}+1}{x}\right)=\prod_{u=0}^{n-1}\left(x^{q^n+1}+x^{q^n}-\left(\beta^{q^{n+u}}+\beta^{q^u}+1 \right)x+1 \right).
$$
Lemma \ref{lemma3N} completes the proof: The polynomial $x^n \psi \left( \displaystyle  \frac{x^{q^n+1}+x^{q^n}+1}{x} \right)$ is irreducible over $\mathbb{F}_q$, since the polynomial $x^{q^n+1}+x^{q^n}-(\beta^{q^{n}}+\beta +1)x+1$ is irreducible over $\mathbb{F}_{q^n}$ and $\deg_q(\beta^{q^n}+\beta+1)=n$.  
\end{pf}

Preliminary versions of Theorems \ref{theorem5},\ref{theorem8} are given in \cite{kyureg-evoyan}.
   
\vspace*{0.3cm}
Further we  use the following result by Sidelnikov to describe  two more composition constructions of explicit families of  irreducible polynomials of degree $n(q^n-1)$  from a given  primitive polynomial of degree $n$.  
\begin{thm1}[Sidelnikov \cite{6Sidelnikov}]
\label{theorem9}
The polynomial 
$$
f(x)=\frac{x^{q+1}-\omega x^q-(x_0+x_1-\omega)x+x_0x_1}{x^2-(x_0+x_1)x+x_0x_1},
$$
\noindent  where $\omega,  x_1, x_0 \in \mathbb{F}_q$, \ $x_0 \not = x_1$,  is irreducible  if and only if  $\displaystyle ~\frac{\omega+x_0}{\omega+x_1}$ is a primitive element of  $\mathbb{F}_q$.
Moreover $f(x)$ has linearly independent roots   over $\mathbb{F}_q$ if  $\omega\not =0$.
\end{thm1}

\begin{thm1}
\label{theorem10}
Let $f(x) \not = x-1$  be a primitive polynomial  of degree $n$ over $\mathbb{F}_{q}$.
Then the polynomial 
$$
F(x)=f\left( x^{q^n}+ x^{q^n-1}\right) \big(f( x+ 1) \big)^{-1}
$$
\noindent of degree $n(q^n-1)$ is irreducible over $\mathbb{F}_{q}$.
\end{thm1}

\begin{pf}
Let $\alpha$ be a  zero of $f(x)$. Then  $\alpha$ is a primitive element of $\mathbb{F}_{q^n}$,
since $f(x)$ is a primitive polynomial of degree $n$ over  $\mathbb{F}_{q}$.
Take  $w=0$, $x_0=\alpha$ and $x_1=1$. 
Note that $x_0= \alpha \not = x_1 =1$ and $\displaystyle \frac{\omega +x_0} {\omega+x_1}= \alpha$  is a primitive element of $\mathbb{F}_{q^n}$. Hence by Theorem \ref{theorem9} 
  the polynomial  
\begin{eqnarray*}
h(x)& = & \frac{x^{q^n+1}- (\alpha+1)x+ \alpha}{x^2-(\alpha+1)x + \alpha} =  \frac{x(x-1)^{q^n}- \alpha (x-1)}{x (x-1)- \alpha (x-1)}\\
 & = & \frac{x(x-1)^{q^n-1}-\alpha}{x-\alpha}
\end{eqnarray*}
is irreducible over $\mathbb{F}_{q^n}$. 
Substituting $x+1$ for $x$  we obtain the polynomial 
$$
g(x) = h(x+1) =\frac {(x+1) x^{q^n-1} -\alpha}{x+(1 -\alpha)}
$$
which is also irreducible over $\mathbb{F}_{q^n}$.  It is easy to see that
$$
(x+1) x^{q^n-1} - \alpha = \big(x-(\alpha-1)\big) \left( x^{q^n-1}+ \alpha x^{q^n-2}+ \cdots + \displaystyle \frac{\alpha}{\alpha -1} \right),
$$
and in particular 
$$
g(x) = x^{q^n-1}+ \alpha x^{q^n-2}+ \cdots + \displaystyle \frac{\alpha}{\alpha -1}.
$$
Since $\alpha$ is a proper element of $\mathbb{F}_{q^n}$ over $\mathbb{F}_{q}$, the degree of the set of coefficients of $g(x)$ over $\mathbb{F}_{q}$ is $n$.
Our next goal is to show that
$$
F(x) = \prod_{u=0}^{n-1} \left ( \frac{(x+ 1)x^{q^n-1} - \alpha ^{q^u}}{x+1- \alpha^{q^u}} \right) = \prod_{u=0}^{n-1} g^{(u)}(x).
$$
Indeed, 
\begin{equation}
\label{eq8a}
f(x) =\displaystyle  \prod _{u=0}^{n-1}(x-\alpha^{q^u})
\end{equation}
over  $\mathbb{F}_{q^n}$. 
Substituting $ (x+1)x^{q^n-1}$, resp. $x+ 1$, for $x$  in (\ref{eq8a}), we obtain  
$$
f \left( (x+1)x^{q^n-1}\right)=\prod_{u=0}^{n-1} \left( (x+1)x^{q^n-1}- \alpha ^{q^u}\right)
$$
and 
$$
f \left(x+1\right)=\prod_{u=0}^{n-1} \left (x+ 1 - \alpha ^{q^u}\right),
$$
which yield
\begin{equation*}
\label{eqN18}
F(x)=\left( f  \left( x+1 \right) \right)^{-1} f \left( (x+1)x^{q^{n-1}} \right)=
\prod_{u=0}^{n-1} \left ( \frac{(x+ 1)x^{q^n-1} -\alpha ^{q^u}}{x+1- \alpha^{q^u}} \right).
\end{equation*}
Finally, the irreducibility of  $F(x)$ over $\mathbb{F}_{q}$ follows from  Lemma \ref{lemma3N}.
\end{pf}

\begin{thm1}
\label{theorem11}
Let $f(x) \not = x-1$ be a primitive polynomial of degree  \ $n$ over $\mathbb{F}_{q}$. Then the polynomial
$$
F(x)= \left( x^{q^n}-2x-1\right)^n f \left( \frac{x^{q^n+1}-x^{q^n}+2x}{x^{q^n}-2x-1}\right)
 \big( (- (x+1))^{n} f(-x) \big)^{-1}
$$
\noindent of degree $n(q^n-1)$ is irreducible over $\mathbb{F}_{q}$.

\end{thm1}

\begin{pf}
 Let $\alpha$ be a zero of $f(x)$.  Thus if   $x_1=-\alpha$, $x_0=-1$ and $\omega =\alpha+1$, then $x_0=-1 \not =x_1=-\alpha$ and $\displaystyle \frac{\omega +x_0}{\omega+x_1}=\frac{\alpha+1-1}{\alpha+1-\alpha}=\alpha$ is a primitive element of $\mathbb{F}_{q^n}$.
 Hence by Theorem \ref {theorem9} the polynomial
 \begin{equation}
\label{eq9}
h(x) = \frac{x^{q^n+1}-x^{q^n}+2x -\alpha (x^{q^n}-2x-1)}{(x+1)(x+\alpha)}
\end{equation}
is irreducible over $\mathbb{F}_{q^n}$. 
Note that 
$$
h(x) = \frac{x^{q^n+1}- (\alpha + 1) x^{q^n}+2x  + 2\alpha x + \alpha }{x ^2 +(\alpha + 1)x + \alpha} = x^{q^n-1} - 2(\alpha + 1) x^{q^n-2} + \ldots + 1,
$$
implying that the degree of the set of coefficients of  $h(x)$ over $\mathbb{F}_{q}$ is equal to $n$ since $\deg_q( -2(\alpha+1)) =n$.

Next we show that $F(x) =  \prod_{u=0}^{n-1} h^{(u)}(x)$ and hence the proof follows from Lemma   \ref{lemma3N}. 
From the  irreducibility of  $f(x)$ over $\mathbb{F}_q$, we have the relation
\begin{equation}
\label{eq10}
f(x)=\prod_{u=0}^{n-1} \left (x-\alpha^{q^u} \right)
\end{equation}  
over $\mathbb{F}_{q^n}$.                                                                              
Substituting $\displaystyle \frac{x^{q^n+1}-x^{q^n}+2x}{x^{q^n}-2x-1}$ for $x$ in (\ref{eq10}) and multiplying both sides of the equation by $\left( x^{q^n}-2x-1\right)^n$, we get
$$
\left(\!x^{q^n}\!-\! 2x \!-\!1 \! \right)^n f \left(\!\displaystyle \frac{x^{q^n+1}\!-\! x^{q^n}+2x}{x^{q^n}\!-\! 2x \!-\!1}\!\right)=\prod_{u=0}^{n-1}\! \left( x^{q^n+1}\!-\! x^{q^n}\! +\! 2x- \! \alpha^{q^u} \left( x^{q^n}\! -\! 2x \!- \! 1 \!\right) \!\right).
$$    
Next substituting $-x$ for $x$ in (\ref{eq10}) and multiplying both sides of the equation by $(-(x+1))^n$, we obtain 
$$
(-(x+1))^n f(-x) =\prod_{u=0}^{n-1}(x+1) (x+\alpha ^{q^n}).
$$
Finally, dividing the first  equation by the second one, we obtain                        
\begin{eqnarray*}
\label{eqN22}
F(x) &=& \left(\frac{\! x^{q^n}\!-\! 2x \! -\! 1\! }{-(x \!+\! 1)}\right)^n  f^{-1}(-x)\, f \! \left(\! \frac{x^{q^n \!+\!1}\!-\! x^{q^n}+2x}{x^{q^n}-\! 2x \! -\! 1}\right) \!\! \nonumber   \\ [4mm]
&=& \prod_{u=0}^{n-1} \! \left( \! \frac{ x^{q^n+1}\! -\! x^{q^n}\! +\! 2x \! - \! \alpha^{q^u} \left(\! x^{q^n}\! -\! 2x \!-\! 1\right)} {(x \!+\! 1)(x \!+\! \alpha^{q^u})}\! \right) = \prod_{u=0}^{n-1} h^{(u)}(x).
\end{eqnarray*}
\end{pf}

\end{document}